\documentclass[11pt]{amsart}
\usepackage{fontawesome5}

\usepackage{amsmath,amsthm,amssymb}
\usepackage{xcolor}
\usepackage{hyperref}

\hypersetup{
colorlinks=false,
linkbordercolor=red,
pdfborderstyle={/S/U/W 1}
}



\usepackage{graphicx}
\usepackage{color}
\usepackage[noend]{algpseudocode}
\usepackage{tikz}

\theoremstyle{plain}
\newtheorem{thm}{Theorem}[section]

\newtheorem{lemma}[thm]{Lemma}

\newcommand{\auto}{\tag{\addtocounter{equation}{1}\theequation\;{\tiny\faCar}}}

\theoremstyle{definition}

\theoremstyle{remark}
\newtheorem{rem}[thm]{Remark}

\numberwithin{equation}{section}

\title[Sharp isoperimetric inequalities]{Sharp isoperimetric inequalities on the Hamming cube II: The critical exponent}

\author[P. Durcik]{Polona Durcik}
\address{Chapman University, Orange, CA, USA}
\email{durcik@chapman.edu}
\author[P. Ivanisvili]{Paata Ivanisvili}
\address{University of California Irvine, Irvine, CA, USA}
\email{pivanisv@uci.edu}
\author[J. Roos]{Joris Roos}
\address{University of Massachusetts Lowell, Lowell, MA, USA}
\email{joris\_roos@uml.edu}
\author[X. Xie]{Xinyuan Xie}
\address{University of California Irvine, Irvine, CA, USA}
\email{xinyuax7@uci.edu}

\date{\today}
\begin{document}
\begin{abstract}
A sharp isoperimetric inequality for the Hamming cube is proved at the critical exponent $\beta=\frac12$. This follows up on previous work, where such bounds were established for $\beta$ near $\frac12$.
As a consequence, this result settles a conjecture of Kahn and Park on cube partitions and yields a sharp $L^1$ Poincar\'e inequality for Boolean-valued functions. It also confirms a low-noise limit for balanced functions predicted by the Hellinger conjecture on noisy Boolean channels in information theory.
\end{abstract}


\maketitle

\section{Introduction}
Let $n\ge 1$ be an integer.
For $A\subset \{0,1\}^n$ and $x\in A$ let
$h_A(x)$ denote the number of edges connecting $x$ with $\{0,1\}^n\setminus A$. In other words, $h_A(x)$ is the number of single-bit flips of $x$ that leave $A$.
If $x\not\in A$, then let $h_A(x)=0$.
Here $\mathbf{E}f=2^{-n} \sum_{x\in \{0,1\}^n} f(x)$ and $|A|=\mathbf{E} \mathbf{1}_A$.

\begin{thm}\label{thm:mainisoperim}
For all $A\subset \{0,1\}^n$ with $|A|\le \frac12$,
\begin{equation}\label{eqn:isoperimhalf}
\mathbf{E} \sqrt{h_A} \ge |A| \sqrt{\log_2(1/|A|)}.
\end{equation}
\end{thm}

This inequality is sharp in two different ways: first, it is an equality when $A$ is a subcube.
Second, there is no corresponding dimension-free lower bound when $\sqrt{h_A}$ is replaced by $h_A^\beta$ for any $\beta<\frac12$, which can be seen by Hamming ball examples.
The study of $\mathbf{E} \sqrt{h_A}$ was first initiated by Talagrand \cite{Tal93} who proved lower bounds such as $\mathbf{E}\sqrt{h_A}\ge \sqrt{2}|A|(1-|A|)$. This was later improved by Bobkov--G\"otze \cite{BG99}, as a consequence of more general results.

Sharp lower bounds for $\mathbf{E} h_A^\beta$ for all $\beta\ge \frac12$ follow from \eqref{eqn:isoperimhalf} and H\"older's inequality; see e.g. \cite[Lemma 5.1]{DIR24}.
Such bounds were previously proved for $\beta\ge 0.50057$ in \cite{DIR24} and prior to that for $\beta\ge \log_2(3/2)\approx 0.585$ in \cite{KP20,BIM23}.
For $\beta=1$ one has Harper's classical isoperimetric inequality \cite{Har66,Ber67,Hart76}.
By the Cauchy-Schwarz inequality, \eqref{eqn:isoperimhalf} implies the following sharp strengthening of the classical isoperimetric inequality for $|A|\le \frac12$,
\begin{equation}\label{eqn:classicalsharpening}
\mathbf{E} h_A \ge \frac{|A|}{|\partial A|} |A|\log_2(1/|A|),
\end{equation}
where $\partial A$ is the support of $h_A$, i.e. the vertex boundary of $A$.
This is an equality if $A$ is a subcube. It improves on the classical isoperimetric inequality because of the factor $\frac{|A|}{|\partial A|}\ge 1$. Moreover, it fails when the factor is replaced by $(\frac{|A|}{|\partial A|})^{\gamma}$ for any $\gamma>1$, again by Hamming ball examples.
Such an inequality was observed by Beltran--Ivanisvili--Madrid \cite[Cor. 1.2]{BIM23} for $\gamma=\frac1{\log_2 (3/2)}-1\approx 0.709$.

The proof of Theorem \ref{thm:mainisoperim} relies on a crucial induction argument of Kahn and Park \cite{KP20}, which reduced proving lower bounds for $\mathbf{E} h_A^\beta$ to finding a certain Bellman function and verifying a two-point inequality.
Our Bellman function is a further refinement of the function used in \cite{DIR24}, which in turn builds on important work by Beltran--Ivanisvili--Madrid \cite{BIM23}.
In the proof of the two-point inequality we draw from estimates and tools introduced in \cite{DIR24}.

\subsection*{Applications}
As an immediate consequence of Theorem \ref{thm:mainisoperim} we settle a conjecture of Kahn and Park \cite[Conjecture 1.3]{KP20} on cube partitions:
\begin{thm}\label{cor:kahnpark}
Let $(A,B,W)$ be a partition of $\{0,1\}^n$ and assume $|A|=\frac12$. Then
\[ |\nabla(A,B)| + \sqrt{n}\,|W|\ge \tfrac12. \]
\end{thm}
Here $|W|=\mathbf{E} \mathbf{1}_W$ and the edge boundary measure $|\nabla(A,B)|$ is defined by $2^{-n} \# \{(x,y)\in \mathcal{E}\,:\,x\in A, y\in B\}$, where $\mathcal{E}$ is the set of edges in the Hamming cube $\{0,1\}^n$.
Theorem \ref{cor:kahnpark} follows from Theorem \ref{thm:mainisoperim} because for $|B\cup W|=|A|=\frac12$, we have
\[\tfrac12\le \mathbf{E} \sqrt{h_{B\cup W}} \le \mathbf{E} (h_{B\cup W}\mathbf{1}_B) + \mathbf{E} (\sqrt{h_{B\cup W}}\mathbf{1}_W),\]
which is no greater than $|\nabla(A,B)| + \sqrt{n} |W|$ as required.

Our next application is a sharp $L^1$ Poincar\'e inequality for Boolean-valued functions.
\begin{thm}\label{thm:poincare} For all $f:\{0,1\}^n\to \{0,1\}$,
\begin{equation}\label{eqn:poincareL1}
\|\nabla f\|_1\ge \|f-\mathbf{E}f\|_1.
\end{equation}
Equality is achieved when $f=\mathbf{1}_A$ for a half-cube $A$.
\end{thm}
Here $|\nabla f(x)|=\sqrt{\sum_{i=1}^n (\tfrac12(f(x)-f(x\oplus e_i)))^2}$ and $\|f\|_1=\mathbf{E} |f|$. Note that it is still open to determine the sharp constant in this inequality for {\it real-valued} $f$, but it is known that this constant lies in the interval $(\frac2{\pi},\sqrt{\frac{2}{\pi}}]$; see \cite{BELP08}, \cite{ILvHV}, also \cite{haonan1}, \cite{Esken1} for quantum and vector-valued inequalities. In particular, \eqref{eqn:poincareL1} shows that the lower bound improves for Boolean-valued functions, which was first observed in \cite{BIM23}.
Theorem \ref{thm:poincare} does not follow directly from Theorem \ref{thm:mainisoperim}, but from the more technical isoperimetric inequality \eqref{eqn:isoperim2} proved below.

Finally, isoperimetric inequalities for the Hamming cube are also related to a family of conjectures in information theory concerning noisy channels (Chen--Nair \cite{CN24}). One of these is the most informative Boolean function conjecture, also known as Courtade--Kumar conjecture \cite{KC13}.
The latter would be implied by a stronger conjecture, known as Hellinger conjecture \cite{ABCJN17, CGN25}. Roughly speaking, these conjectures concern optimality questions for Boolean functions under standard noisy-channel models on $\{0,1\}^n$, where each bit is independently perturbed.
In the low-noise limit, the Hellinger conjecture would in particular imply that
\begin{equation}\label{eqn:sensitivity}
\mathbf{E}\sqrt{s_f} \ge 1
\end{equation}
holds for balanced Boolean functions $f:\{0,1\}^n\to \{0,1\}$ (i.e. $\mathbf{E}f=\tfrac12$), where $s_f(x)$ is the sensitivity of $f$ at $x$, which is defined as the number of single-bit flips of $x$ that change the value of $f(x)$, i.e. $s_f(x)=h_A(x)+h_{A^c}(x)$.
This implication is well-known, but for convenience of the reader we provide details in \S \ref{sec:it}.
Since $h_A$ and $h_{A^c}$ have disjoint support, $\sqrt{s_f}=\sqrt{h_A}+\sqrt{h_{A^c}}$, so Theorem \ref{thm:mainisoperim} with $|A|=\frac12$ immediately implies \eqref{eqn:sensitivity}.

\subsection*{Acknowledgments} P.D., P.I. and J.R. thank the American Institute of Mathematics (AIM) for funding our SQuaRE project. The authors were also supported in part by grants from the National Science Foundation CAREER-DMS-2152401 (P.I.), DMS-2154835 (J.R.) and by the Simons Foundation
SFI-MPS-TSM-00013943 (P.D.), SFI-MPS-TSM-00014075 (J.R.), Simons Fellowship, and Humboldt Research Fellowship (P.I.).

\section{Preliminaries}\label{sec:envelope}

\subsection{Reduction to a two-point inequality}

For $0\le x\le y\le 1$ and $B:[0,1]\to [0,\infty)$ let
\begin{equation}\label{eqn:Gdef}
G_1[B](x,y) = \sqrt{(y-x)^{2}+B(y)^{2}} + B(x) - 2 B(\tfrac{x+y}2),
\end{equation}
\[ G_2[B](x,y) = y-x + (\sqrt{2}-1) B(y) + B(x) - 2 B(\tfrac{x+y}2),\]
and $G = \max(G_1, G_2).$
Kahn and Park \cite{KP20} proved that if $B$ satisfies $B(0)=B(1)=0$ and the two-point inequality
\[G[B](x,y)\ge 0\]
holds for all $0\le x\le y\le 1$, then $\mathbf{E} \sqrt{h_A}\ge B(|A|)$
for all $A\subset \{0,1\}^n$ and $n\ge 1$.

This reduces the proof of Theorem \ref{thm:mainisoperim} to finding a Bellman function $B$ and verifying the two-point inequality.
We use the following refinement of the Bellman function used in \cite{DIR24}.
For $x\in (0,1]$ define
\[ L(x) = x \sqrt{\log_2(1/x)} \]
and $L(0)=0$.
Let $Q(x)$ be the unique cubic interpolating polynomial such that $Q(0)=Q(1)=0$, $Q(\frac12)=\frac12$ and $Q(\frac14)=2^{-3/2}$, i.e.
\[ Q(x) = \tfrac23 x (1-x)(2^{5/2}-3+4(3-2^{3/2}) x). \]
The idea of using the polynomial $Q$ goes back to Beltran--Ivanisvili--Madrid \cite{BIM23}.
Let $I(x)$ denote the Gaussian isoperimetric profile, i.e. the unique function on $[0,1]$ such that
$I(0)=I(1)=0, I\cdot I''=-1$.
For a parameter $w\in (\frac12,1]$ let
\[ {\mathrm{J}}_w(x) = \tfrac12 I(\tfrac1{2w})^{-1} I(\tfrac{1-x}{w}). \]
Note that ${\mathrm{J}}_w(\frac12)=\frac12$ and ${\mathrm{J}}_w(1)=0$.
Then define
\begin{equation}\label{eqn:bbdef}
B_w(x) = \left\{ \begin{array}{ll}
L(x) & \text{for } x\in [0, \frac14],\\
Q(x)
& \text{for }x\in [\frac14, \frac12],\\
{\mathrm{J}}_w(x) & \text{for }x\in [\frac12, 1].
\end{array}\right.
\end{equation}
For some motivation on the choice of such functions, see \cite[\S 2]{DIR24}.
In the following we shall fix
\begin{equation}\label{eqn:w1def}
w=w_1=\tfrac{29}{32}
\end{equation}
and write ${\mathrm{J}}={\mathrm{J}}_{w_1}$.
Note that since $B_w(x)\ge L(x)$ for $x\in [0,\frac12]$ (see \cite[Lemma 5.9 (2)]{DIR24}), in order to prove Theorem \ref{thm:mainisoperim} it suffices to prove that there exists $w$ such that
\begin{equation}\label{eqn:maintwopt}
G[B_w](x,y)\ge 0
\end{equation}
for all $0\le x\le y\le 1$.
As a consequence of this argument we obtain the isoperimetric inequality
\begin{equation}\label{eqn:isoperim2}
\mathbf{E} \sqrt{h_A} \ge B_{w_1}(|A|),
\end{equation}
valid for all $A\subset \{0,1\}^n$.

\begin{rem}The result remains valid for an interval of values $w$.
It could be interesting to further investigate the optimal (smallest) value $w_*$ of $w$, but we will not pursue this here. Numerical evidence suggests that $w_*\in [0.897, 0.898]$.
\end{rem}

\begin{rem}
It is also interesting to ask what the sharp behavior of $\mathbf{E}\sqrt{h_A}$ is as $|A|\to 1-$.
While our new isoperimetric inequality \eqref{eqn:isoperim2} improves on the corresponding inequality proved in \cite[(2.7)]{DIR24} for all values of $|A|\le \frac12$, this changes when $|A|>\frac12$.
Specifically, for $|A|$ near $1$, the inequality \cite[(2.7)]{DIR24} is slightly stronger than \eqref{eqn:isoperim2}.
Similarly, for all values of $|A|>\frac12$, the inequality \cite[(2.6)]{DIR24}, which gives sharp lower bounds for $\mathbf{E}h_A^\beta$ when $\beta\ge \beta_0=0.50057$, is stronger than the estimate $\mathbf{E} h_A^\beta\ge {\mathrm{J}}(|A|)$ that follows from \eqref{eqn:isoperim2}. This is consistent with the fact that ``Case $J$'' in \cite[\S 6]{DIR24} is more difficult to handle than Case $J$ in the present paper (see \S \ref{sec:caseJ}). For $\beta\ge\beta_0$ and $|A|>\frac12$, the best known lower bound is that given in \cite[Thm. 1.6]{DIR24}, which further improves on \cite[(2.6)]{DIR24} as $|A|\to 1-$.
\end{rem}

\subsection{\texorpdfstring{Properties of the function ${\mathrm{J}}$}{Properties of the function J}}\label{sec:prelimJ}
We record a few properties of the function ${\mathrm{J}}$ for later use.
The function ${\mathrm{J}}={\mathrm{J}}_{w_1}$ satisfies
\[{\mathrm{J}}\cdot {\mathrm{J}}'' = -\gamma\quad \text{with}\;\gamma=(4w_1^2I(\tfrac{1}{2w_1})^2)^{-1}.\]
In particular, ${\mathrm{J}}$ is concave.
Note that $\gamma\in [1.945, 1.946]$.
Let
\[x_1 = 1-\tfrac{w_1}{2}= \tfrac{35}{64} = 0.546875.\]
For $x\in (0,\tfrac{1}{64}]$ we have
\begin{equation}\label{eqn:Jlower}
\mathrm{J}(1-x)>x\sqrt{\log(w_1/x)}.
\end{equation}

Indeed, the Gaussian isoperimetric profile satisfies the lower bound
\[I(x)\ge \sqrt{2}\cdot x\sqrt{\log(1/x)} \Big(1-\tfrac{1}{2}\tfrac{\log \log (1/x)}{\log (1/x)} - \tfrac{\log (2\sqrt{\pi})}{\log(1/x)} \Big), \]
valid for $x\in (0,\frac15]$ (see \cite[Proposition 5.11]{DIR24} for the proof).
Thus,
\[{\mathrm{J}}(1-x) \ge \tfrac{\sqrt{2}}{2w_1}I(\tfrac{1}{2w_1})^{-1} (1-\varepsilon) x\sqrt{\log(w_1/x)},\]
where
\[\varepsilon = \tfrac{1}{2}\tfrac{\log \log (w_1/x)}{\log (w_1/x)} + \tfrac{\log (2\sqrt{\pi})}{\log(w_{1}/x)}. \]
Then \eqref{eqn:Jlower} follows by recognizing that $\varepsilon$ is an increasing function of $x$ and evaluating at $x=\frac1{64}$.

The function ${\mathrm{J}}$ is increasing on $[\tfrac{1}{2},x_1]$ and decreasing on $[x_1,1]$. Therefore, an interval enclosure for ${\mathrm{J}}$ is given by
\begin{equation}
\label{eqn:J_intencl}
\underline{{\mathrm{J}}}(\underline{x},\overline{x}) = \min({\mathrm{J}}(\underline{x}),{\mathrm{J}}(\overline{x})),
\end{equation}
\[\overline{{\mathrm{J}}}(\underline{x},\overline{x}) = \left\{ \begin{array}{ll}
{\mathrm{J}}(\overline{x}),& \textup{if } \overline{x}<x_1 \\
{\mathrm{J}}(\underline{x}), & \textup{if } \underline{x}>x_1\\
{\mathrm{J}}(x_1),&\textup{else.}
\end{array}\right.\]
Also observe ${\mathrm{J}}'(x_1)=0$ and ${\mathrm{J}}''<0$, so ${\mathrm{J}}'$ is strictly decreasing and $|{\mathrm{J}}'|$ has the interval enclosure
\[ \underline{|{\mathrm{J}}'|}(\underline{x},\overline{x})= \left\{
\begin{array}{ll}
{\mathrm{J}}'(\overline{x}), & \text{if}\;\overline{x}<x_1,\\
-{\mathrm{J}}'(\underline{x}), & \text{if}\;\underline{x}>x_1,\\
0, & \text{else,}
\end{array}
\right.\]
\[ \overline{|{\mathrm{J}}'|}(\underline{x},\overline{x}) = \max(|{\mathrm{J}}'(\underline{x})|, |{\mathrm{J}}'(\overline{x})|). \]

\subsection{\texorpdfstring{Properties of the functions $L$ and $Q$}{Properties of the functions L and Q}}
We will use the following straightforward facts about the functions $L$ and $Q$. See \cite[Lemma 5.6, Lemma 5.8]{DIR24} for details on the proofs.
\begin{lemma}[L]\label{lem:L}
The function $x\mapsto L(x)$ is strictly concave on $[0,1]$. If $x\in (0,e^{-\sqrt{3}/2})$, then $L'''(x)>0$.
\end{lemma}

\begin{lemma}[Q]
\label{lem:Q}
If $x\in [0,\tfrac12]$, then
$Q(x)\ge 0$, $Q'(x)>0$, and $Q''(x)<0$.
\end{lemma}

\subsection{Computer-assisted proofs}
For the proof of the two-point inequality \eqref{eqn:maintwopt}
we make use of computer-assisted verification using interval arithmetic, more specifically using recursive dyadic partitioning as explained in \cite[\S 4]{DIR24}. The verification code is implemented using \emph{FLINT/Arb}, an open source library for arbitrary precision interval arithmetic \cite{Arb}, \cite{Flint}.
It can be accessed on GitHub at
\[
\texttt{\href{https://github.com/roos-j/dir24-isoperim}{https://github.com/roos-j/dir24-isoperim}}
\]
Claims verified in this manner will be tagged with ({\tiny\faCar}). Single point evaluations are not included in this treatment. We stress that it is straightforward to verify all these claims independently (see \cite[\S 4]{DIR24} for more details; see also \cite{Tuc11}, \cite{Rump}, \cite{GomezSerrano}).

\section{Proof of the two-point inequality}\label{sec:maintwopt}
In this section we prove \eqref{eqn:maintwopt}. The argument is split into six cases according to the piecewise definition of $B_w$.
We follow some of the arguments in \cite[\S 6]{DIR24} with the most important differences in the critical case, Case $J$,
and in Cases $LJ$, $QJQ$.

\subsection{\texorpdfstring{Case $J$: $\frac12\le x\le y\le 1$}{Case J}}
\label{sec:caseJ}
This is the main case, where the failure for $\beta$ near $\frac12$ occurred in \cite{DIR24}.
It suffices to show
\begin{equation}
\label{eq:caseJ}
\sqrt{(y-x)^2 + {\mathrm{J}}(y)^2} + {\mathrm{J}}(x)\ge 2 {\mathrm{J}}( \tfrac{x+y}{2}).
\end{equation}

We give an argument
inspired by Bobkov's two-point inequality \cite{Bob97} (see also \cite[Prop. 3.3]{DIR24}).
Set $c=\tfrac{x+y}{2}$ and $u=\tfrac{y-x}{2}$. Then
it suffices to show
\[ \sqrt{4u^2 + {\mathrm{J}}(c+u)^2} \ge 2 {\mathrm{J}}( c) - {\mathrm{J}}(c-u)\]
for all $0\le u\le \min(c-\frac12,1-c)=c_*$ and $c \in [\frac12, 1]$.
By concavity of ${\mathrm{J}}$, the right-hand side is non-negative so we may equivalently show
\[4u^2 + {\mathrm{J}}(c+u)^2 \ge 4 {\mathrm{J}}(c)^2 + {\mathrm{J}}(c-u)^2 -4 {\mathrm{J}}(c){\mathrm{J}}(c-u).\]
Letting
\[F_c(u) = {\mathrm{J}}(c+u)^2 - {\mathrm{J}}(c-u)^2 + 4u^2 - 4 {\mathrm{J}}( c)^2 +4 {\mathrm{J}}(c){\mathrm{J}}(c-u), \]
we compute
\[F_c'(u) = 2 {\mathrm{J}}(c+u){\mathrm{J}}'(c+u) + 2{\mathrm{J}}(c-u){\mathrm{J}}'(c-u) +8u -4{\mathrm{J}}(c){\mathrm{J}}'(c-u),\]
\[F_c''(u) = 2 ({\mathrm{J}}'(c+u)^2-{\mathrm{J}}'(c-u)^2) + 8 - 4\gamma\tfrac{{\mathrm{J}}(c)}{{\mathrm{J}}(c-u)},\]
where we used that ${\mathrm{J}}\cdot {\mathrm{J}}'' = -\gamma$.
It suffices to show that
\begin{equation}\label{eqn:FcDDpos}
F_c''(u) \ge 0
\end{equation}
for $u\in [0,c_*]$ and $c\in [\frac12, 1]$.
Indeed, since $F_c(0)=F'_c(0)=0$ this implies $F_c(u)\ge 0$ for $u\in [0,c_*]$ as required.
To show \eqref{eqn:FcDDpos} we distinguish two cases.

\subsubsection*{\texorpdfstring{Case I: $x=c-u\ge x_1$.}{Case I}}
Since ${\mathrm{J}}$ is decreasing on $[x_1, 1]$ and $\gamma\le 2$,
\[ F''_c(u) \ge 2 ({\mathrm{J}}'(y)^2-{\mathrm{J}}'(x)^2) .\]
Since $({\mathrm{J}}')^2$ is convex (see \cite[Lemma 3.2]{DIR24}),
\[ {\mathrm{J}}'(y)^2 - {\mathrm{J}}'(x)^2 \ge 2u (({\mathrm{J}}')^2)'(x). \]
Observe $(({\mathrm{J}}')^2)' = 2 {\mathrm{J}}' \cdot {\mathrm{J}}'' = - 2 \gamma {\mathrm{J}}' / {\mathrm{J}} \ge 0$ on the interval $[x_1,1]$, since ${\mathrm{J}}'\le 0$ there. This shows that $F_c''(u)\ge 0$.

\subsubsection*{\texorpdfstring{Case II: $x=c-u\in [\frac12, x_1]$}{Case II}}
We use monotonicity properties of the functions ${\mathrm{J}},|{\mathrm{J}}'|$ on the intervals $[\frac12, x_1]$ and $[x_1,1]$ (see \S \ref{sec:prelimJ}) to estimate
\[ F''_c(u) = 2{\mathrm{J}}'(y)^2 - 2{\mathrm{J}}'(x)^2 + 8 - 4 \gamma {\mathrm{J}}((x+y)/2)/{\mathrm{J}}(x) \]
\[ \ge -2{\mathrm{J}}'(\tfrac12)^2 + 8 - 4\gamma {\mathrm{J}}(x_1)/{\mathrm{J}}(\tfrac12)\in [0.08, 0.09]. \]

\subsection{\texorpdfstring{Case $Q$: $0\le x\leq y \leq\frac12$}{Case Q}}
\label{sec:caseQ}
In \cite[Prop. 6.2]{DIR24} it was proved that $G_{1}[\max(L,Q)]\ge 0$ in this region, which implies the claim.

\subsection{\texorpdfstring{Case $LJQ$: $0\leq x \leq \tfrac{1}{4},\,\tfrac{1}{2}\leq y \leq 1,x+y\leq 1$.}{Case LJQ}}
\label{sec:caseLJQ}
We split the region into two pieces and treat each separately.\\

\noindent {\it Case I:} $\tfrac{1}{16}\leq x\leq \tfrac{1}{4}, \tfrac{1}{2}\leq y\leq \tfrac34$.
\label{sec:LJQ-I}
In this rectangle, we use the tight lower bound for $G[B_w](x,y)$ given by
\[ \underline{h_{LJQ,1}}(\underline{x},\overline{x},\underline{y},\overline{y})=
\max\Big(\sqrt{(\underline{y}-\overline{x})^{2}+\underline{{\mathrm{J}}}(\underline{y},\overline{y})^{2}}, \underline{y}-\overline{x}+(\sqrt{2}-1) \underline{{\mathrm{J}}}(\underline{y},\overline{y})\Big) \]
\[ + L(\underline{x}) - 2 Q(\tfrac{\overline{x}+\overline{y}}{2}).
\]
By dyadic partitioning we obtain the lower bound
\begin{equation}\label{eqn:h_LJQ_1_auto}\auto
G[B_w](x,y)>10^{-6}
\end{equation}
for $(x,y)\in [\frac1{16},\frac14]\times[\frac12,\frac34]$. \\

\noindent {\it Case II.}
\label{sec:LJQ-II}
For the remainder of the case $LJQ$ it suffices to show that if
$x\in [0,\frac14],$ $y\in [\frac12,1]$, and
if in addition $x\le \frac1{16}$ or $y\ge \frac34$ holds, then
\[ G_{LJQ}(x,y)=y-x + (\sqrt{2}-1) {\mathrm{J}}(y) +L(x)-2Q(\tfrac{x+y}{2})\geq 0.\]

For fixed
$y\in [0, 1]$, the function $x\mapsto G_{LJQ}(x,y)$ is strictly concave on $[0,\frac14]$.
Indeed, in \cite[Lemma 6.7]{DIR24} it was shown that
\begin{equation}\label{eqn:LJQlempf1}
\partial_x^2 G_{LJQ}(x,y) = L''(x) - \tfrac12 Q ''(\tfrac{x+y}2) <0.
\end{equation}
Hence it remains to verify the non-negativity of $G_{LJQ}$ on the $x$-boundary of the region. First we verify
$G_{LJQ}(0,y)\ge 0$ for all $y\in [\frac12,1]$. Let
\[ f(y)=G_{LJQ}(0,y).
\]
Using ${\mathrm{J}}\cdot{\mathrm{J}}''=-\gamma$, we compute
\[ f'(y) = 1 + (\sqrt{2}-1){\mathrm{J}}'(y) - Q'(\tfrac{y}2), \]
\[f''(y) = -\gamma (\sqrt{2}-1) {\mathrm{J}}(y)^{-1} -\tfrac12 Q''(\tfrac{y}2).\]
Using
\begin{equation}
\label{J_der4}
{\mathrm{J}}^{(4)} = -\gamma(\gamma+2({\mathrm{J}}')^2){\mathrm{J}}^{-3},
\end{equation}
we compute $f^{(4)}(y)<0$.
Hence $f''$ is strictly concave on $(\frac12,1)$. Since
\begin{equation}
f''(\tfrac{1}{2})>0.03,\quad f''(\tfrac{9}{16})>0.05,\quad
f''(\tfrac{4}{5})<-0.2,
\end{equation}
the function $f''$ has exactly one zero $y_0$ on $(\tfrac{1}{2},1)$, where $f'$
achieves its unique maximum.
Since $f'$ is increasing on $(\tfrac{1}{2},y_0)$ and decreasing on $(y_0,1)$, and
\[f'(\tfrac{9}{16})>0.05,\quad
f'(\tfrac{15}{16})<-0.08,\]
it follows that $f'$ has exactly one zero $y_1$ on $(\tfrac{1}{2},1)$. Therefore, $f$ is increasing on $[\frac12,y_1]$ and decreasing on $[y_1,1]$. Together with $f(1)=f(\tfrac{1}{2})=0$, this implies $f(y)\geq 0$ for $y\in [\frac12,1]$.

Next we show $G_{LJQ}(\tfrac1{16},y)\ge 0$ for all $y\in[\tfrac12,\tfrac34]$. A tight lower bound is
\[ \underline{h_{LJQ,2}}(\underline{y},\overline{y}) = \underline{y}-\tfrac{1}{16} + (\sqrt{2}-1) \underline{{\mathrm{J}}}(\underline{y},\overline{y}) +L(\tfrac{1}{16})-2Q(\tfrac{1}{32}+ \tfrac{\overline{y}}{2}).\]
Dyadic partitioning then yields
\begin{equation}\label{eqn:h_LJQ_2_auto}\auto
G_{LJQ}(\tfrac1{16},y)>10^{-4}
\end{equation}
for all $y\in [\frac12,\frac34]$.

Finally, we show $G_{LJQ}(1-y,y)\ge 0$ for all $y\in [\tfrac34,1]$.
The function $y\mapsto G_{LJQ}(1-y,y)$ is concave on $[\tfrac34,1]$ since
\[\tfrac{d^2}{dy^2}G_{LJQ}(1-y,y)= (\sqrt{2}-1) {\mathrm{J}}''(y)+L''(1-y)\le 0.\]
Here we used ${\mathrm{J}}''\leq 0$ and $L''\le 0$ by Lemma \ref{lem:L}.
Thus, it suffices to show non-negativity of the function at $y=1$, and $y=\tfrac{3}{4}$, in which cases
\[G_{LJQ}(0,1)=0,\quad G_{LJQ}(\tfrac{1}{4},\tfrac{3}{4})> 0.01.\]

\subsection{\texorpdfstring{Case $LJ$: $0 \leq x \leq \tfrac{1}{4}, \tfrac{1}{2}\leq y\leq 1, 1\leq x+y$}{Case LJ}}
\label{sec:caseLJ}
Expressed in the variables $u=(y-x)/2$ and $c=(x+y)/2$ it suffices to show
\[ h_{LJ}(u,c) = 2u + (\sqrt{2}-1){\mathrm{J}}(c+u)+L(c-u)-2{\mathrm{J}}(c)\ge 0 \]
for all $c\in [\frac12, \frac58]$ and $u\in [c-\frac14,1-c]$.
Using ${\mathrm{J}}\cdot {\mathrm{J}}''=-\gamma$ we compute
\[ \partial_u^2 h_{LJ}(u,c) = -\gamma(\sqrt{2}-1) {\mathrm{J}}(c+u)^{-1} + L''(c-u)\le 0, \]
since $L$ is concave by Lemma \ref{lem:L}.
Therefore it suffices to evaluate $h_{LJ}(u,c)$ at the boundaries $u=c-\frac14$ and $u=1-c$.\\

\noindent {\it Case I:} $u=c-\frac14$. In this case it suffices to show that
\[ h_{LJ,1}(c) = h_{LJ}(c-\tfrac14,c)>0 \]
for all $c\in [\frac12,\frac58]$.
A tight lower bound for $h_{LJ,1}$ is given by
\[ \underline{h_{LJ,1}}(\underline c,\overline c) = 2\underline{c} - \tfrac12 + (\sqrt{2}-1) \underline{{\mathrm{J}}}(2\underline{c}-\tfrac14,2\overline{c}-\tfrac14) + L(\tfrac14) - 2 \overline{{\mathrm{J}}}(\underline c,\overline c).\]
Dyadic partitioning yields
\begin{equation}\label{eqn:h_LJ_1_auto}\auto
h_{LJ,1}(c) > 0.01\quad\text{for all}\; c\in [\tfrac12,\tfrac58].
\end{equation}

\noindent {\it Case II:} $u=1-c$. In this case it suffices to show that
\[ g(c) = 2(1-c) + L(2c-1) - 2{\mathrm{J}}(c) \ge 0 \]
for all $c\in [\frac12,\frac58]$.
The proof is complete if we show that $g$ is increasing on this interval since $g(\frac12)=0$.
To this end we compute
\[ g'(c) = -2 + 2 L'(2c-1) - 2{\mathrm{J}}'(c). \]
We check that $g'(c)>0$ for $c\in (\frac12,\frac58]$. Note that $L'$ and ${\mathrm{J}}'$ are both decreasing.
Thus for $c\in(\frac12,\frac{9}{16}]$ we have
\[ g'(c) \ge -2 + 2 L'(\tfrac18) - 2 {\mathrm{J}}'(\tfrac12) \ge 0.26 \]
and for $c\in [\frac{9}{16}, \frac{19}{32}]$ we have
\[ g'(c) \ge -2 + 2 L'(\tfrac3{16}) - 2 {\mathrm{J}}'(\tfrac{9}{16}) \ge 0.3 .\]
Finally, for $c\in [\frac{19}{32},\frac58]$,
\[ g'(c) \ge -2 + 2 L'(\tfrac14) - 2 {\mathrm{J}}'(\tfrac{19}{32}) \ge 0.17. \]
Thus $g'(c)\ge 0.17$ for all $c\in(\frac12,\frac58]$.

\subsection{\texorpdfstring{Case $QJQ$: $\frac14\leq x\leq \frac12\leq y, x+y\leq 1$}{Case QJQ}}
\label{sec:caseQJQ}
We will show that for all $(x,y)$ in the larger region $[\frac14,\frac12]\times [\frac12,\frac34]$,
\begin{equation}\label{eqn:QJQclaim}
(y-x)^2 + {\mathrm{J}}(y)^2 - (2Q(\tfrac{x+y}{2})-Q(x))^2 \geq 0.
\end{equation}

The $y$-derivative of the left-hand side in \eqref{eqn:QJQclaim} is two times
\begin{equation}\label{eqn:gQJQdef}
h_{QJQ,1}(x,y)=y-x +{\mathrm{J}}(y){\mathrm{J}}'(y)-(2Q(\tfrac{x+y}{2})-Q(x))Q'(\tfrac{x+y}{2}).
\end{equation}

We begin by showing that this quantity is strictly positive for all $(x,y)\in [\frac14,\frac12]\times [\frac{1}{2},\frac{33}{64}]$, thus reducing to the case $y=\frac12$ on this region.
This can be done by dyadic partitioning. In order to formulate a tight lower bound we record the monotonicity of the various terms appearing in \eqref{eqn:gQJQdef}.

First, one observes that the quantity $2Q(\frac{x+y}2)-Q(x)$ is positive, increasing in $y$ and decreasing in $x$ (for details see \cite[\S 6.5]{DIR24}). Also, $Q'$ is decreasing and positive by Lemma \ref{lem:Q}. Note positivity continues to hold on the interval $[0,\tfrac{33}{64}]$ as one verifies by a single point evaluation.
Finally, the function $x\mapsto {\mathrm{J}}(x){\mathrm{J}}'(x)$ is decreasing on $x\in [\frac12, \frac34]$. Indeed, \[({\mathrm{J}}\cdot {\mathrm{J}}')' = ({\mathrm{J}}')^2-\gamma.\]
Since $({\mathrm{J}}')^2$ is convex (see \cite[Lemma 3.2]{DIR24}), it suffices to evaluate $({\mathrm{J}}')^2-\gamma$ at the endpoints $x=\frac12$ and $x=\frac34$ which shows $({\mathrm{J}}')^2-\gamma<-1<0$.

Therefore, a tight lower bound of $h_{QJQ,1}$ is given by
\[ \underline{h_{QJQ,1}}(\underline{x},\overline{x},\underline{y},\overline{y}) = \underline{y}-\overline{x} + {\mathrm{J}}(\overline{y}){\mathrm{J}}'(\overline{y}) - (2Q(\tfrac{\underline{x}+\overline{y}}2)-Q(\underline{x}))Q'(\tfrac{\underline{x}+\underline{y}}2). \]
Dyadic partitioning shows that
\begin{equation}\label{eqn:h_QJQ_1_auto}\auto h_{QJQ,1}>10^{-5}
\end{equation}
on $[\frac14,\frac12]\times [\frac{1}{2},\frac{33}{64}]$.

To finish the proof for $y\le \frac{33}{64}$, it now suffices to show the inequality \eqref{eqn:QJQclaim} for $y=\frac12$, that is,
\[(\tfrac{1}{2}-x)^2 + \tfrac{1}{4} - (2Q(\tfrac{x}{2} + \tfrac{1}{4})-Q(x))^2 \ge 0 \]
for $x\in [\frac14,\frac12]$.
This inequality was shown in \cite[\S 6.5]{DIR24}.

It remains to consider the case $(x,y)\in [\frac{1}{4},\frac{1}{2}]\times [\tfrac{33}{64},\tfrac{3}{4}]$. A tight lower bound of the left-hand side in \eqref{eqn:QJQclaim} is
\[\underline{h_{QJQ,2}}(\underline x,\overline x, \underline y, \overline y)=(\underline{y}-\overline{x})^2 + \underline{{\mathrm{J}}}(\underline y, \overline{y})^2 - (2Q(\tfrac{\underline{x}+\overline{y}}{2})-Q(\underline{x}))^2.\]
By dyadic partitioning we obtain
\begin{equation}\label{eqn:h_QJQ_2_auto}\auto \text{LHS of }\eqref{eqn:QJQclaim}>10^{-7}
\end{equation}
on $[\frac14,\frac12]\times [\tfrac{33}{64},\tfrac{3}{4}]$.

\subsection{\texorpdfstring{Case $QJ$: $\frac14\le x\le \frac12\le y\le 1$, $x+y\ge 1$}{Case QJ}}
\label{sec:caseQJ}
We split into two subcases.

\subsubsection{Near diagonal: $3/8\leq x\leq 1/2, 1/2\leq y\leq 5/8, x+y\geq 1$}
\label{sec:QJ-I}
Here we show that for $\frac38\le x\le \frac12\le y\le \frac58$,
$G_1[B_w](x,y)\ge 0$.
This will follow from positivity of an equivalent expression
\[ g_{QJ}(x,y)=(y-x)^{2} + {\mathrm{J}}(y)^2- (2 {\mathrm{J}}(\tfrac{x+y}{2})-Q(x))^2. \]

Note
that
$g_{QJ}(\tfrac12, y)\ge 0$
for all $y\in [\tfrac12,\tfrac58]$.
Indeed, this is equivalent to showing $G_1[{\mathrm{J}}](\frac12,y)\ge 0$ for these values, which already follows from Section \ref{sec:caseJ}.
Thus, to finish the proof of this case it remains to show that $\partial_x g_{QJ}$ is negative, or equivalently, the positivity of
\[ -\tfrac{1}{2}\partial_x g_{QJ}(x,y) = y-x+(2{\mathrm{J}}(\tfrac{x+y}{2})-Q(x)) ({\mathrm{J}}'(\tfrac{x+y}{2})-Q'(x)).\]

Recalling Lemma \ref{lem:Q} and the enclosures for ${\mathrm{J}}, |{\mathrm{J}}'|$ (see \eqref{eqn:J_intencl} and the subsequent displays), a tight lower bound of this expression can be given by
\[ \underline{h_{QJ,1}}(\underline{x},\overline{x},\underline{y},\overline{y})
=\underline{y}-\overline{x}+(2\underline{{\mathrm{J}}}(\tfrac{\underline{x}+\underline{y}}2,\tfrac{\overline{x}+\overline{y}}2)-Q(\overline{x})) {\mathrm{J}}'(\tfrac{\overline{x}+\overline{y}}2)\mathbf{1}_{(\overline{x}+\overline{y})/2<x_1}
\]
\[ - (2\overline{{\mathrm{J}}}(\tfrac{\underline{x}+\underline{y}}2,\tfrac{\overline{x}+\overline{y}}2)-Q(\underline{x})) \overline{|{\mathrm{J}}'|}(\tfrac{\underline{x}+\underline{y}}2,\tfrac{\overline{x}+\overline{y}}2)\mathbf{1}_{\text{not}\;(\overline{x}+\overline{y})/2<x_1}
- (2\overline{{\mathrm{J}}}(\tfrac{\underline{x}+\underline{y}}2,\tfrac{\overline{x}+\overline{y}}2)-Q(\underline{x}))Q'(\underline{x}). \]
Dyadic partitioning yields
\begin{equation}\label{eqn:h_QJ_1_auto}\auto
-\tfrac{1}{2}\partial_x g_{QJ} > 10^{-5}\quad\text{on}\quad[\tfrac14,\tfrac12]\times [\tfrac12,\tfrac58].
\end{equation}

\subsubsection{Far from diagonal: $\tfrac14\leq x\leq \tfrac12, \frac58\leq y\le 1$, $x+y\ge 1$}
\label{sec:QJ-II}
It suffices to show that for all $x\in [\frac14, \frac12]$ and $y\in[\frac58,1]$,
\begin{equation}\label{eqn:h_QJ_2_auto}\auto
h_{QJ,2}(x,y)= \sqrt{(y-x)^2 +{\mathrm{J}}(y)^2}+Q(x)-2{\mathrm{J}}(\tfrac{x+y}{2})>10^{-7}.
\end{equation}
This follows by dyadic partitioning using the tight lower bound
\[ \underline{h_{QJ,2}}(\underline{x},\overline{x},\underline{y},\overline{y}) = \sqrt{(\underline{y}-\overline{x})^2 + {\mathrm{J}}(\overline{y})^2} + Q(\underline{x})-2\overline{{\mathrm{J}}}(\tfrac{\underline{x}+\underline{y}}{2},\tfrac{\overline{x}+\overline{y}}{2}). \]

\section{Proof of the Poincar\'e inequality}\label{sec:poincare}
In this section we prove Theorem \ref{thm:poincare} using \eqref{eqn:isoperim2}.
Every Boolean-valued function $f$ can be written as $f=\mathbf{1}_A$ for some $A\subset \{0,1\}^n$.
Then $\|f-\mathbf{E} f\|_1 = 2|A| (1-|A|)$ and
$\|\nabla f\|_1 = \tfrac{1}{2} ( \mathbf{E} \sqrt{h_A} + \mathbf{E} \sqrt{h_{A^c}})$.
By \eqref{eqn:isoperim2} it therefore suffices to show that
\begin{equation}\label{eqn:poincarepenult}
G_P(x)= \tfrac{1}{2}( B_{w_1}(x) + B_{w_1}(1-x)) - 2x(1-x) \ge 0.
\end{equation}
for all $x\in [0,1]$.
Since $G_P(x)=G_P(1-x)$ it suffices to show this for $x\in [0,\tfrac12]$.

\subsubsection*{Case I: $x\in [0,\tfrac1{64}]$}
We must show
\[\tfrac{1}{2}( L (x) + {\mathrm{J}}(1-x)) - 2x (1-x) \ge 0.\]
By \eqref{eqn:Jlower}, the left-hand side is
\[ \ge x (\tfrac{1}{2} \sqrt{\log_2(1/x)} + \tfrac{1}{2} \sqrt{\log(w_1/x)}- 2) = x\cdot g_{P,1}(x). \]
Since $g_{P,1}$ is a decreasing function and
$g_{P,1}(\tfrac1{64})>0.2$
we obtain $G_{P}(x)\ge 0.2\, x$ for all $x\in [0,\tfrac1{64}]$.

\subsubsection*{Case II: $x\in [\tfrac1{64},\tfrac14]$}
Since $x_1<\tfrac34$ the function $x\mapsto {\mathrm{J}}(1-x)$ is increasing on $[\frac1{64},\frac14]$.
Therefore a tight lower bound for $G_{P}(x)$ is given by
\[\underline{h_{P,1}}(\underline{x},\overline{x})= \tfrac{1}{2}(L(\underline{x}) + {\mathrm{J}}(1-\underline{x})) - 2\overline{x}(1-\underline{x}).\]
Dyadic partitioning shows that
\begin{equation}\label{eqn:h_P_1_auto}\auto
G_{P}(x)>10^{-4}
\end{equation}
for all $x\in [\frac1{64},\frac14]$.

\subsubsection*{Case III: $x\in [\tfrac14,\tfrac12]$}
Here we show that
\[G_P(x) = \tfrac{1}{2}( Q(x) + {\mathrm{J}}(1-x)) - 2x(1-x) \ge 0.\]
The left-hand side vanishes at $x=\tfrac12$. Thus, it suffices to show that the function on the left is decreasing in $x$, i.e. that
\[h_{P,2}(x)=-\partial_x G_{P}(x)>0\]
for all $x\in [\frac14,\frac12]$.
We compute
\[ h_{P,2}(x)=-\tfrac{1}{2} Q'(x)+\tfrac{1}{2} {\mathrm{J}}'(1-x) - 4x+2. \]
The function $x\mapsto {\mathrm{J}}'(1-x)$ is increasing on $[0, \tfrac12]$
and $Q'$ is decreasing, so
a tight lower bound for $h_{P,2}$ is given by
\[ \underline{h_{P,2}}(\underline{x},\overline{x}) = - \tfrac{1}{2} Q'(\underline{x}) + \tfrac{1}{2} {\mathrm{J}}'(1-\underline{x}) -4\overline{x}+2 .\]
By dyadic partitioning,
\begin{equation}\label{eqn:h_P_2_auto}\auto
h_{P,2}(x) > 10^{-4}
\end{equation}
for all $x\in[\frac14,\frac12]$, which finishes the proof of \eqref{eqn:poincarepenult}.

\section{Connection to the Hellinger conjecture}\label{sec:it}
We briefly explain how the Hellinger conjecture would also imply \eqref{eqn:sensitivity}.

Let $f:\{0,1\}^n\to \{\pm 1\}$ be balanced, i.e.\ $\mathbf{E} f =0$.
In this case the Hellinger conjecture \cite{ABCJN17,CN24,CGN25} asserts that for every correlation parameter $\rho \in [-1,1]$,
\begin{equation}\label{eqn:hellinger}
\mathbf{E}\sqrt{1 - (T_\rho f)^2} \ge \sqrt{1-\rho^2},
\end{equation}
where $T_\rho f$ is the standard noise operator (see \cite{ODonnell}). Let $p=\tfrac{1-\rho}{2}\in [0,1]$ and let $Z=(Z_1,\dots,Z_n)$ be a random vector in $\{0,1\}^n$ with i.i.d. coordinates satisfying $\mathbf{P}(Z_j=1)=p$, so that
\[ T_\rho f(x) = \mathbf{E}_{Z}\,f(x\oplus Z).\]
We are interested in the low-noise regime as $p\to 0+$.
Conditioning on the event of at most a single bit flip gives
\[ T_\rho f(x) = (1-p)^n f(x) + p (1-p)^{n-1} \sum_{j=1}^n f(x\oplus e_j) + R_p(x), \]
where $R_p(x)=\mathbf{E}_Z [f(x\oplus Z)\mathbf{1}_{|Z|\ge 2}]$ (here $|Z|=\sum_{j=1}^n Z_j$).
Since $f$ is Boolean valued,
\[ |R_{p}(x)| \le \mathbf{P}(|Z|\ge 2) \le \binom{n}{2} p^2. \]
Also, $\sum_{j=1}^n f(x\oplus e_j) = (n-2s_f(x))f(x)$,
where $s_f(x)$ is the sensitivity of $f$ at $x$, i.e. the number of indices $j$ so that $f(x\oplus e_j)\not=f(x)$.
Thus,
\[
T_\rho f(x)
= f(x)(1-2p\, s_f(x))+O(n^2 p^2),
\]
by expanding the binomial $(1-p)^n$ and using the bound on $|R_{p}(x)|$.
This implies in particular that
\[ \sqrt{1-(T_\rho f(x))^2} = 2\sqrt{p(1-p)} \sqrt{s_f(x)} + O(n p)\]
for every $x\in \{0,1\}^n$.
Dividing by $\sqrt{1-\rho^2}=2\sqrt{p(1-p)}$ and taking expectations over $x$ we see that
\eqref{eqn:hellinger} would imply
\[ \mathbf{E} \sqrt{s_f} \ge 1 - O(n \sqrt{p}) \]
so that for every fixed $n$ we obtain \eqref{eqn:sensitivity} in the limit $p\to 0+$.

\end{document}